\documentclass[11pt]{article}
\usepackage{latexsym}
\usepackage{graphicx}
\usepackage{amsmath}
\usepackage{amssymb}
\usepackage{color}
\usepackage{epsfig}
\usepackage{cite}

\newtheorem{problem}{Problem}[section]

\newtheorem{lemma}{Lemma}
\newtheorem{theorem}[problem]{Theorem}
\newtheorem{prop}[problem]{Proposition}

\newtheorem{example}{Example}[section]

\def\ref{\noindent{\bf Reference}}

\makeatletter      
\@addtoreset{equation}{section}
\makeatother       

\title{On the number of certain Galois extensions of  local fields }
\author{Da-sheng Wei$^1$ and Chun-gang Ji$^2$ \\
\small{$^1$Department of Mathematics, the University of Science
and Technology of China, Hefei 230026.} \\
\small{{\it Email:} dshwei@ustc.edu}\\
\small{$^2$Department of Mathematics, Nanjing Normal University, Nanjing 210097.}\\
\small{{\it Email:} cgji@njnu.edu.cn}}
\date{}

\begin{document}%
\maketitle \noindent{\bf Abstract:} In this paper, we will
calculate the number of Galois extensions of local fields with
Galois group $A_n$ and $S_n$.



\section{Introduction}
Let $p$ be a prime, $F$ a finite extension of $p$-adic field
$\mathbb{Q}_p$ with $[F:\mathbb{Q}_p]=m$. Let $k$ be its residue
field with $[k:\mathbb{Z}/p\mathbb{Z}]=f$. It's clear $f|m$. Let
$\pi$ be a uniformizer of $F$ and $e$ be the absolute ramification
index of $F/Q_p$. Hence $m=ef$. And let $\mu_l$ denote the set of
$l$-th roots of unity. All notations are standard if not
explained.

Since the number of the extensions of the local field with the
given degree is finite, see [2]. One can ask for the formula that
gives the number of extensions of a given degree. Krasner [1] gave
such a formula, and Serre [5] also computed the number of
extensions using a different method in the proof of his famous
``mass formula". Paul and Roblot [6] gave another proof for the
formula. Similarly one can also ask for the formula that gives the
number of Galois extensions of the given degree. In particular, it
is possible to ask for the formula that gives the number of the
Galois extensions with prescribed finite Galois group $G$. We
denote the number by $\nu(F,G)$. If $G$ is a $p$-group with $\mu_p
\nsubseteq F$, $\check{\text{S}}\text{afarevi}\check{\text{c}}$
[3] gave an explicit formula for the number of the $G$-extensions
of $F$:
$$\nu(F,G)=\frac{1}{|Aut(G)|}(\frac{|G|}{p^d})^{m+1}\prod_{i=1}^{d-1}(p^{m+1}-p^i),$$
where $d$ is the minimal number of generators of $G$. If $G$ is a
$p$-group, and $\mu_p \subset F$, Yamagishi [7] obtained a formula
for $\nu(F,G)$.

In this paper, we will calculate the number of $S_n$-extensions
and $A_n$-extensions of $F$, where $S_n$ is $n$-th symmetric group
an $A_n$ is $n$-th alternating group, that is,

\begin{theorem}
$(1)$ Suppose the prime $p\neq 3$, then
$$\nu(F,S_3)=\begin{cases} 0 \ \ \   &\text{if}\ \ \mu_3 \subset F, \\
 1\ \ \  & \text{if}\ \ \mu_3 \not \subset F.\end{cases}$$
$(2)$ Suppose $p=3$, then
$$\nu(F,S_3)=\begin{cases} 3^{m+1}-3 \ \ \  & \text{if}\ \ \mu_3 \subset F,\\
3^m+\frac{3^{m+1}}{2}-\frac{3}{2} \ \ \   &\text{if}\ \ \mu_3 \not
\subset F.\end{cases}$$
\end{theorem}

\begin{theorem}
$(1)$ Suppose the prime $p \geqslant 3$, then
$$\nu(F,S_4)=\nu(F,A_4)=0.$$
$(2)$ Suppose $p=2$, then
$$\nu(F,A_4)=\begin{cases} 4(2^{2m}-1)/3 \ \ \   &\text{if}\ \mu_3 \subset F, \\
 (2^{2m}-1)/3 \ \ \  & \text{if}\ \mu_3 \not
\subset F;\end{cases}$$
$$\nu(F,S_4)=\begin{cases} 0 \ \ \   &\text{if}\ \mu_3 \subset F, \\
 2^{2m+1}-1 \ \ \  & \text{if}\ \mu_4 \nsubseteq F\
 \text{and}\ m\ \text{is\ even}\ \text{and}\ f=1, \\
 2^{2m}-1\ \ \  & \text{otherwise}.\end{cases}$$
\end{theorem}

\section{Some lemmas}
\begin{lemma}
Suppose $n\geqslant 5$, there does not exist a Galois extension of
$F$ with Galois group $S_n$ or $A_n$.
\end{lemma}
\begin{pf} If $F$ is a local field, then the Galois group $\text{Gal}(K/F)$ is a solvable group
for every finite Galois extension $K$ of $F$. It is well-known
that $A_n$ and $S_n$ is not a solvable group for $n\geqslant
5$.$\Box$
\end{pf}

It is well-known that the number of $S_2$-extensions and
$A_3$-extensions of local field by local classfield throry. So we
only need to calculate the number of Galois extensions of field
$F$ with Galois group $S_3$ ,$S_4$ and $A_4$. The following lemma
plays important role in our calculation.
\begin{lemma}
Let $K$ be a Galois extension of $F$ with the Galois group $G$.
For any subgroup $A$ of $G$, let $F_A$ be the field fixed by $A$,
then the Galois closure $cl(F_A)$ of $F_A$ is a subfield of $K$
and $\text{Gal}(K/cl(F_A))=N_{G}(A)$, where $N_{G}(A)=\{s\in
A|t^{-1}st\in A, \text{for\ all}\ t \in G\}$.
\end{lemma}

We can get some Galois extensions from  some non-Galois extensions
by its' Galois closure. For example, if G=$S_3$, $D_8$, $A_4$ and
$S_4$, we can choose  $A$ to be order $2$, order $2$, order $3$
and order $6$ non-normal subgroup respectively of $G$, where $D_8$
is the $2$-sylow subgroup of $S_4$. We know $A$ is isomorphic to
$\mathbb{Z}/2\mathbb{Z}$, $\mathbb{Z}/2\mathbb{Z}$,
$\mathbb{Z}/3\mathbb{Z}$ and $S_3$ respectively. And their
normalization is trivial. By the above lemma, the Galois
extensions of $F$ with Galois group $S_3$ can be gotten by the
Galois closure of extension of degree $3$ of $F$, and the Galois
extensions of $F$ with Galois group $D_8$, $A_4$ and $S_4$ can be
gotten by the Galois closure of extensions of degree $4$ of $F$.

Let $M(n)$ be the set of all extensions of degree $n$ of $F$. Let
$Ab(n)$ be the set of abelian extensions of degree $n$ of $F$. And
let $M(G)$ be the set of Galois extensions of $F$ with the Galois
group $G$. Let $K$ be the Galois closure of an extension of degree
$n$ of $F$ . The Galois group $\text{Gal}(K/F)$ is a subgroup of
$S_n$. Obviously the order of $Gal(K/L)$ must be divided by $n$.
So there are the following two maps:
$$f: \ \ M(3)\rightarrow Ab(3)\cup M(S_3)$$ and
$$g: \ \ M(4)\rightarrow Ab(4)\cup M(D_8)\cup M(A_4)\cup M(S_4)$$ by
$$L\rightarrow cl(L).$$
 The two maps are surjective. Any inverse image $L$ of an
element $K$ of $M(G)$ is a subfield of $K$ and $L$ is not a Galois
extension of $F$ if $G$ is not an abelian group. By Galois theory,
the Galois group $\text{Gal}(L/K)$ is not a normal subgroup of
$G$. For $G=S_3, D_8,A_4$ and $S_4$, we consider respectively the
number of the non-normal subgroup isomorphic to
$\mathbb{Z}/2\mathbb{Z}$,$\mathbb{Z}/2\mathbb{Z}$,$\mathbb{Z}/3\mathbb{Z}$
and $S_3$. The subgroups are non-normal except there is an
order-$2$ normal subgroup in $D_8$. So the number of the inverse
image of any element of $M(S_3), M(D_8), M(A_4)$ and $M(S_4)$ are
$3,4,4,4$ respectively. Let $|S|$ denote cardinality of a finite
set $S$. Let $\nu(F,G)$ denote the number of $M(G)$. So there is
the following result.
\begin{lemma}
$$|M(3)|=|Ab(3)|+3\nu(F,S_3),$$
$$|M(4)|=|Ab(4)|+4\nu(F,D_8)+4\nu(F,A_4)+4\nu(F,S_4).$$
\end{lemma}

\section{The proof of theorem 1.1}
In the following, we denote $m=[F:\mathbb{Q}_p]$ and $e$ is the
absolute ramification index of $F$ and $q=p^f$ is the number of
elements of the residue field of $F$.

\begin{pf}(1) For $p\neq 3$,
$$M(3)=\{K|\ [K:F]=3,K\ \text{is\ a\ tamely\ ramified\ extension\ of\ } F\}.$$
(i) If $3$-th roots of unity are contained in $F$, then
$M(3)=Ab(3)$. By lemma 3 $$\nu(F,S_3)=0.$$ (ii)If $3$-th roots of
unity are not contained in $F$, then $|M(3)|=4$, $|Ab(3)|=1$. By
lemma 3
$$\nu(F,S_3)=1.$$
(2) For $p=3$, by Krasner'theorem [3],
$$|M(3)|=3q^e+6(q-1)(\sum_{a=0}^{e-1}
q^a)+1=9q^e-5.$$
Suppose $\mu_3 \not \subset F$, then
$$|Ab(3)|=\frac{1}{2}(\frac{3}{3})^{m+1}(3^{m+1}-1)=\frac{3^{m+1}-1}{2}=\frac{3q^e-1}{2}.$$
Suppose $\mu_3 \subset F$, then
$$|Ab(3)|=4.$$
By lemma 3
$$\nu(F,S_3)=\begin{cases} \frac{5q^e-3}{2} \ \ \   &\text{if}\ \mu_3 \not \subset F, \\
3q^e-3 \ \ \  & \text{if}\ \mu_3 \subset F .\end{cases}$$ $\Box$
\end{pf}

\section{The proof of theorem 1.2}
First we give some propositions.
\begin{prop}
If the prime $p \geqslant 3$, then
$$\nu(F,S_4)=\nu(F,A_4)=0.$$
\end{prop}
\begin{pf} Suppose $K$ is a Galois extension of $F$ with Galois group
$S_4$. There must exist  intermediate fields $F^{tr}$ and $F^{ur}$
such that $\text{Gal}(K/F^{tr})$ is a $p$-group,
$\text{Gal}(F^{tr}/F^{ur})$ and $\text{Gal}(F^{ur}/F)$ are cyclic
groups. By Galois theory, there is a $p$-group $S'$ which is a
normal subgroup of $S_4$. Since $p\geqslant 3$, $S'$ must be
$(1)$. Since $S_4$ has not a cyclic normal subgroup $S$ such that
$S_4/S$ is also cyclic, a contrary is derived. The result follows
similarly for $A_4$.$\Box$
\end{pf}

\begin{prop} If $p=2$,
$$\nu(F,A_4)=\begin{cases} 4(2^{2m}-1)/3 \ \ \   &\text{if}\ \mu_3 \subset F,\\
 (2^{2m}-1)/3 \ \ \  & \text{if}\ \mu_3 \not \subset F  .\end{cases}$$
\end{prop}
\begin{pf}Let $K$ be an $A_4$-extension of $F$. Since $K_4$ is a
normal subgroup of $A_4$, there exists a (unique) Galois subfield
$F'$ of degree $3$ of $K$, where $K_4 \cong
\mathbb{Z}/2\mathbb{Z}\times \mathbb{Z}/2\mathbb{Z}$. By [2],
$$|F'^*/(F'^*)^2|=4q^{3e},$$
$$|F^*/(F^*)^2|=4q^e.$$
It is clear that the natural map of $F^*/(F^*)^2\rightarrow
F'^*/(F'^*)^2$ is an injection. We consider the action on
$F'^{*}/(F'^{*})^2$ of the Galois group $\text{Gal}(F'/F)$ .\\

$(1)$ Denote $G'=Gal(F'/F)$. Then the following result holds.
$$(F'^*/(F'^*)^2)^{G'}
 = F^*/(F^*)^2.$$  In the following we will prove it. Let $a \in
F'^*/(F'^*)^2$ be a fixed point of $\text{Gal}(F'/F)$ and $a\neq
1$ in $F'^*/(F'^*)^2$. Then $F'(\sqrt{a'})$ is a Galois extension
of $F$, where $a'$ represents a lifting of $a$ in $F'^*$. There
isn't an order-$2$ normal subgroup in $S_3$, so
$$\text{Gal}(F'(\sqrt{a'})/F)=\mathbb{Z}/6\mathbb{Z}.$$ Let $F''$ be
the fixed field of the normal subgroup $\mathbb{Z}/3\mathbb{Z}$.
There exists an element $b\in F^*/(F^*)^2$, such that
$$F''=F(\sqrt{b}).$$ Then $$F'(\sqrt{a'})=F'(\sqrt{b}).$$ So $a=b$ in
$F'^*/(F'^*)^2$.\\
\\
(2) Let $\sigma$ be a generator of $G'$. If $x
\notin F^*/(F^*)^2$. Then there are the following two case:\\
(i) $N_{F'/F}(x) =1\ \text{in}\ F^*/(F^*)^2$,\\
(ii) $N_{F'/F}(x)\neq 1\ \text{in}\ F^*/(F^*)^2$.\\
In (i), the field $F'(\sqrt{x},\sqrt{\sigma x})$ is an
$A_4$-extension of $F$ since the Galois group of
$F'(\sqrt{x},\sqrt{\sigma x})/F$ isomorphic to $K_4\rtimes
\mathbb{Z}/3\mathbb{Z}$, where $K_4 \cong
\mathbb{Z}/2\mathbb{Z}\times \mathbb{Z}/2\mathbb{Z}$.\\
In (ii),
$$F'(\sqrt{x},\sqrt{\sigma
x},\sqrt{\sigma^2 x})=F'(\sqrt{x/\sigma(x)},\sqrt{\sigma
x/\sigma^2 x},\sqrt{N_{F'/F}(x)}),$$ So $F'(\sqrt{x},\sqrt{\sigma
x},\sqrt{\sigma^2 x})$ is a
$(A_4\times\mathbb{Z}/2\mathbb{Z})$-extension of $F$. And for any
$(A_4\times\mathbb{Z}/2\mathbb{Z})$-extension $K$, it is generated
by an $A_4$-extension of $F$ and an  extension of degree $2$ of
$F$. Denote the unique extension of degree 3 by $F'$, then there
exist $x \in F'^*/(F'^*)^2$ and $a\in F^*/(F^*)^2$ satisfying
$x\notin F^*/(F^*)^2$, $N_{F'/F}(x)\in (F^*)^2$ and $a \notin
(F^*)^2$, such that $K=F'(\sqrt{x},\sqrt{\sigma x},\sqrt{a})$. It
is easy to see
$$K=F'(\sqrt{ax/\sigma x},\sqrt{a\sigma x/\sigma^2 x},\sqrt{a\sigma^2
x/x})$$ since $N_{F'/F}(x)=x\sigma x\sigma^2 x=1\ \text{in}\
$$ F^*/(F^*)^2$, then $\sigma x/\sigma^2x=x\sigma^2 x/\sigma^2 x=x\ \text{in}\
$$ F'^*/(F'^*)^2$. \\
Consider $K$ as an extension of $F'$, there exist $7$ subfields
with order $2$ over $F'$ which are one-to-one correspondent to
$\{y,\sigma y,\sigma^2 y,y\sigma y,\sigma y\sigma^2 y,y\sigma^2
y,N_{F'/F}(y) \}$, where $y=ax/\sigma x$. And $N_{F'/F}(y)=a \neq
1\ \text{in}\ F^*/(F^*)^2$. The $\text{Gal}(F'/F)$-orbits are
$\{y,\sigma y,\sigma^2 y\}$, $\{y\sigma y,\sigma y\sigma^2
y,y\sigma^2 y\}$, and $\{N_{F'/F}(y) \}$. So any
$(A_4\times\mathbb{Z}/2\mathbb{Z})$-extension of $F$
can be generated by the element of (ii) and the claim follows.\\
$$\nu(F,A_4\times\mathbb{Z}/2\mathbb{Z})=\nu(F,A_4)\nu(F,\mathbb{Z}/2\mathbb{Z}).$$
(i) Suppose $\mu_3 \not \subset F$, $F'$ is the unique unramified
extension of degree $3$ of $F$, so
$$3\nu(F,A_4)+3\nu(F,A_4\times\mathbb{Z}/2\mathbb{Z})=4q^{3e}-4q^e,$$
$$\nu(F,A_4\times\mathbb{Z}/2\mathbb{Z})=\nu(F,A_4)(4q^e-1).$$ Then
$$\nu(F,A_4)=(q^{2e}-1)/3.$$
(ii) Suppose $\mu_3 \subset F$, $F'$ is the unique  unramified
extension and $3$ totally ramified extensions of degree  $3$ of
$F$, so
$$3\nu(F,A_4)+3\nu(F,A_4\times\mathbb{Z}/2\mathbb{Z})=4(4q^{3e}-4q^e),$$
$$\nu(F,A_4\times\mathbb{Z}/2\mathbb{Z})=\nu(F,A_4)(4q^e-1).$$ Then
$$\nu(F,A_4)=4(q^{2e}-1)/3.$$
$\Box$
\end{pf}

\begin{prop} Let $p=2$,
$$\nu(F,S_4)=\begin{cases} 0 \ \ \   &\text{if}\ \mu_3 \subset F,\\
 2^{2m+1}-1 \ \ \  & \text{if}\ \mu_4 \nsubseteq F\
 \text{and}\ m\ \text{is\ even}\ \text{and}\ f=1,\\
 2^{2m}-1\ \ \  & \text{otherwise}.\end{cases}$$
\end{prop}
\begin{pf} By Krasner's theorem [3],
$$|M(4)|=16q^{3e}-4q^{2e}-5.$$
By the local classfield throry and the dual theory of the finite
abelian group, the following equation holds:
$$|Ab(4)|=|\{S:\ S \text{\ is\ the\ subgroup \ of \ order\ 4\ of\ } F^*/(F^*)^4\}|.$$
Let  $T_1$ be the set consisting of the elements of order
$\leqslant2$ in $F^*/(F^*)^4$, $T_2$ be the set consisting of the
elements of order $4$ in $F^*/(F^*)^4$. The sequence
$$0\rightarrow T_1 \rightarrow F^*/(F^*)^4\rightarrow (F^*)^2/(F^*)^4\rightarrow 0$$
is exact. The third map is $a \mapsto a^2$. So
$$|T_1|=|F^*/(F^*)^4|/|(F^*)^2/(F^*)^4|.$$
Suppose $\mu_4\nsubseteq F$, then
$$|T_1|=4q^e,$$
$$|T_2|=8q^{2e}-4q^e.$$
Then
$$|Ab(4)|=|T_2|/2+(|T_1|-1)(|T_1|-2)/6=20q^{2e}/3-4q^e+1/3.$$
Suppose $\mu_4\subset F$, then
$$|T_1|=4q^e,$$
$$|T_2|=16q^{2e}-4q^e.$$
Then
$$|Ab(4)|=|T_2|/2+(|T_1|-1)(|T_1|-2)/6=32q^{2e}/3-4q^e+1/3.$$
Since $D_8$ is a $2$-group, by theorem 2.2 [7],
$$\nu(F,D_8)=\begin{cases} q^e(q^e-1)(4q^e-1) \ \ \   & \text{if}\ \mu_4\subset F\ \text{or}\ \mu_4 \nsubseteq F\
 \text{and}\ m\ \text{is\ even}\ \text{and}\ f=1,\\
 q^e(2q^e-1)^2 \ \ \  & \text{otherwise}.\end{cases}$$
 By lemma 3,\\
(i) If $\mu_{4} \subset F$, then
$$\nu(F,S_4)=(|M(4)|-|Ab(4)|)/4-\nu(F,A_4)-\nu(F,D_8)=4(q^{2e}-1)/3-\nu(F,A_4).$$
(ii) If $\mu_{4} \nsubseteq F$, and $m$ is odd or $m$ is even and
$f\geqslant 2$, then
$$\nu(F,S_4)=(|M(4)|-|Ab(4)|)/4-\nu(F,A_4)-\nu(F,D_8)=4(q^{2e}-1)/3-\nu(F,A_4).$$
(iii) If $\mu_{4} \nsubseteq F$ and  $m$ is even and $f=1$, then
$$\nu(F,S_4)=(|M(4)|-|Ab(4)|)/4-\nu(F,A_4)-\nu(F,D_8)=(7q^{2e}-4)/3-\nu(F,A_4).$$
By proposition 4.2, so
$$\nu(F,S_4)=\begin{cases} 0 \ \ \   &\text{if}\ \mu_3 \subset F,\\
 2q^{2e}-1 \ \ \  & \text{if}\  \mu_4 \nsubseteq F\
 \text{and}\ n\ \text{is\ even}\ \text{and}\ f=1,\\
 q^{2e}-1\ \ \  & \text{otherwise}.\end{cases}$$
 $\Box$
\end{pf}

Remark: Since $K_4$ is a normal subgroup of $S_4$ and
$S_4/K_4\cong S_3$, there exists an $S_3$-subextension in an
$S_4$-extension of $F$ by Galois theory. If $\mu_3\subset F$ and
$p\neq 3$, then $\nu(F,S_3)=0$. So $\nu(F,S_4)=0$. This gives
another proof for a case of $\nu(F,S_4)$.

Using these propositions, the proof of theorem 1.2 is obtained.
This completes the proof of theorem 1.2.

\section{Examples}
\begin{example} Let $F=\mathbb{Q}_p$, then\\
$(1)$
$$\nu(\mathbb{Q}_p,S_3)=\begin{cases} 6 \ \ \   &\text{if}\ p=3,\\
0 \ \ \   &\text{if}\ p\equiv 1\ \  \mod 3,\\
1 \ \ \  & \text{if}\ p  \equiv 2\ \  \mod 3 .\end{cases}$$
$(2)$
$$\nu(\mathbb{Q}_p,A_4)=\begin{cases} 1 \ \ \   &\text{if}\ p=2,\\
 0 \ \ \  & \text{if}\ p>2  .\end{cases}$$
$(3)$
$$\nu(\mathbb{Q}_p,S_4)=\begin{cases} 3 \ \ \   &\text{if}\ p=2,\\
 0 \ \ \  & \text{if}\ p>2  .\end{cases}$$
$(4)$
$$\nu(\mathbb{Q}_p,S_n)=\nu(\mathbb{Q}_p,A_n)=0\ \ (n\geqslant 5).$$
\end{example}

\end{document}